\documentclass[a4paper,12pt]{article}
\usepackage[english]{babel}
\usepackage[T2A]{fontenc}
\usepackage[cp1251]{inputenc}
\usepackage[english]{babel}
\usepackage{amsthm}
\usepackage[tbtags]{amsmath}
\usepackage{amsfonts,amssymb}
\begin{document}

\textbf{Centrally Essential Semigroup Algebras}

\hfill \textbf{Oleg Lyubimtsev},\footnote{Nizhny Novgorod State University}

\hfill \textbf{Askar Tuganbaev}\footnote{National Research University `MPEI'}

\textbf{Abstract.} For a cancellative semigroup S and a field F, it is proved that the semigroup algebra FS is centrally essential if and only if  the group of fractions $G_S$ of the semigroup $S$  exists and the group algebra $FG_S$ of $G_S$ is centrally essential. The semigroup algebra of a cancellative semigroup is centrally essential if and only if it has the classical right ring of fractions which is a centrally essential ring. There exist non-commutative centrally essential semigroup algebras over fields of zero characteristic (this contrasts with the known fact that  centrally essential group algebras over fields of zero characteristic are commutative).\\ The work of O. Lyubimtsev is done under the state assignment No~0729-2020-0055.



\section{Introduction}\label{sec1}

We consider only associative rings which are not necessarily unital. An associative ring $R$ is said to be \textsf{centrally essential} if either $R$ is commutative or for every non-central element $a$, there exist non-zero central elements $x,y$ with $ax = y$. If the ring $R$ with center $Z(R)$ has the non-zero identity element, then $R$ is centrally essential if and only if the module $R_{Z(R)}$ is an essential extension of the module $Z(R)_{Z(R)}$. It is clear that any commutative ring is centrally essential. Every semiprime or right non-singular centrally essential ring is commutative; see \cite[Theorem 1.5]{MT19b}. In \cite[Proposition 2.4]{MT19b}, it is proved that all idempotents of a centrally essential ring are central. Centrally essential group algebras over fields of characteristic $0$ are commutative; see \cite[Remark 1.2]{MT18}. However, a centrally essential ring can be non-commutative. For example, there exist finite non-commutative centrally essential group algebras over fields of prime characteristic; see \cite{MT18}. In addition, there exist Abelian torsion-free groups such that their endomorphism rings are non-commutative centrally essential rings; see \cite{LT20}. In 
\cite[Theorem 1]{LT21b}, it is proved that a group algebra is a centrally essential if and only if it has the classical right ring of fractions which is a centrally essential ring.

In this paper, we study centrally essential semigroup algebras. The main results of the paper are Theorem 1.1 and Theorem 1.2.

\textbf{Theorem 1.1.}\\
\textbf{1.} Let $S$ be a cancellative semigroup and let $F$ be a field. The semigroup algebra $FS$ is centrally essential if and only if  the group of fractions $G_S$ of the semigroup $S$ exists and the group algebra $FG_S$ of $G_S$ is centrally essential.

\textbf{2.} There exist non-commutative centrally essential semigroup algebras over fields of zero characteristic (this contrasts with the known fact that  centrally essential group algebras over fields of zero characteristic are commutative).

The following theorem extends the result of \cite[Theorem 1]{LT21b} to semigroup algebras of cancellative semigroups.

\textbf{Theorem 1.2.}\\
The semigroup algebra of a cancellative semigroup is centrally essential if and only if it has the classical right ring of fractions which is a centrally essential ring.

In what follows $F$ is a field, $S$ is a semigroup, $FS$ is  the semigroup algebra of the semigroup $S$ over the field $F$. The center of the semigroup $S$ or the semigroup algebra $FS$ are denoted by $Z(S)$ and $Z(FS)$, respectively. If $a = \sum\alpha_ss\in FS$, then 
$\text{supp}(a) = \{s\in S\mid\,\alpha_s\neq 0\}$. For any two elements $a$ and $b$ of a ring, we set $[a, b] = ab - ba$. We give some definitions used in the paper.

A semigroup $S$ is said to be a \textsf{left cancellative semigroup} if for any $a, b, c\in S$, it follows from $ca = cb$ that $a = b$. \textsf{Right cancellative semigroups} are defined dually. A left and right cancellative semigroup is called a \textsf{cancellative semigroup}. It is well known that a torsion cancellative semigroup is a group; e.g., see, \cite{CP61}. A cancellative semigroup can be embedded in the right group of fractions if and only if the intersection of any two principal right ideals of the semigroup $S$ is not empty, i.e., $sS\cap tS\neq \emptyset$ for all $s, t\in S$ (the right Ore condition). If $S$ also satisfies the left Ore condition, which is defined similarly, then the group $G_S = SS^{-1} = S^{-1}S$ is called the \textsf{group of fractions} of the semigroup $S$. Every element of the group $G_S$ can be written in the form $a^{-1}b$ and in the form $cd^{-1}$; $a, b, c, d\in S$. 

If $G$ is a group, then we denote 
$$
\Delta(G) = \{x\in G\mid \,\,\, \mid x^G \mid < \infty\} = (x\in G\mid \,\,\, \mid G:C_G(x)\mid < \infty\}.
$$
$\Delta(G)$ is a characteristic subgroup of the group $G$. In \cite[Lemma 2.1]{MT18}, it is proved that if the group ring $RG$ is centrally essential, then the group $G$ is an $FC$-group, i.e., $G = \Delta(G)$. 

Let $S$ be a cancellative semigroup and $s\in S$. If for some $x\in S$ there exists $t\in S$ such that $xs = tx$, then the element $t$ is uniquely defined; it is denoted by $s^x$. Then $\Delta(S)$ is the set of elements $s\in S$ such that the elements $s^x$ are defined for all $x\in S$ and the set 
$\{s^x\mid \, x\in S\}$ is finite. If $s\in \Delta(S)$, then we assume that $D_S(s) = \{s^x\mid \, x\in S\}$. It is clear that if $S$ can be embedded in the group of fractions $G_S$, then for $s\in \Delta(S)$, the set $D_S(s)$ can be embedded in the set of conjugate elements for $s$ in $G_S$. If $S$ is a cancellative semigroup, then $Z(FS)$ is the $F$-subspace in $FS$ generated by the elements of the form $\sum_{t\in D_S(s)}t$, where 
$s\in \Delta(S)$; see \cite[Theorem 9.10]{Ok91}.

An element $r$ of the ring $R$ is called a \textsf{left non-zero-divisor} or a \textsf{right regular element} if it follows from $rx = 0$ that $x = 0$ for any $x\in R$. In any centrally essential ring, one-sided zero-divisors are two-sided zero-divisors; see \cite[Lemma 2.2]{MT20b}. A ring $R$ has the \textsf{right} (resp., \textsf{left}) \textsf{classical ring of fractions} $Q_{\text{cl}}(R_r)$ (resp., $Q_{\text{cl}}(R_l)$) if and only if for any two elements $a, b\in R$, where $b$ is a non-zero-divisor, there exist  elements $c, d\in R$, where $d$ is a non-zero-divisor, such that $bc = ad$ (resp., 
$cb = da$). If  there exist the rings $Q_{\text{cl}}(R_r)$ and $Q_{\text{cl}}(R_l)$, then they are isomorphic over $R$. In this case, one says that there exists the \textsf{classical ring of fractions} $Q_{\text{cl}}(R)$. It is well known that every commutative ring has the commutative the classical ring of fractions.

\section{Proof of Theorem 1.1}\label{sec2}

\textbf{Proposition 2.1.}\label{prop2.1} 
Let $S$ be a cancellative semigroup. If the semigroup algebra $FS$ is a centrally essential ring, then $S = \Delta(S)$.

$\lhd$ By assumption, if $s\in S$ then we have $0\neq cs = d$ for some $c, d\in Z(FS)$. Then for any $y\in \text{supp}(d)$ there exists 
$x\in \text{supp}(c)$ such that $xs = y$. It follows from \cite[Proposition 9.2(iii)]{Ok91} that $x, y\in \Delta(S)$. In addition, $\Delta(S)$ is a right and left Ore set in $S$ and $G_{\Delta(S)} = \Delta(S)^{-1}\Delta(S) = \Delta(S)\Delta(S)^{-1}$ is an $FC$-group; see \cite[Corollary 9.6 and Proposition 9.8(iii)]{Ok91}. Consequently, $s = x^{-1}y$, where $x^{-1}\in \Delta(S)^{-1}$, $y\in \Delta(S)$. For any $t\in S$, we have $x^t\in \Delta(S)$. Therefore, it follows from $tx = x^tt$ that $(x^t)^{-1}t = tx^{-1}$, i.e., $(x^t)^{-1} = (x^{-1})^t$ in the group $G_{\Delta(S)}$. Then the element 
$s^t = (x^{-1}y)^t = (x^{-1})^ty^t$ there exists for any $t\in S$; see \cite[Basic property (a), P. 108]{Ok91}. Further,
$$
\{s^t\mid \, t\in S\} = \{(x^{-1}y)^t\mid \, t\in S\} = \{(x^{-1})^t\mid \, t\in S\}\cdot \{y^t\, \mid \, t\in S\} = 
$$
$$
\{(x^t)^{-1}\mid \, t\in S\}\cdot \{y^t\, \mid \, t\in S\}.
$$
The first set is finite, since set $\{x^t\mid \, t\in S\}$ is finite. It follows from $y\in \Delta(S)$ that the second set is finite. Therefore, the set $D_S(s)$ is finite and $s\in\Delta(S)$.~$\rhd$

\textbf{Corollary 2.2.}\label{cor2.2}
If $FS$ is a centrally essential semigroup algebra of the cancellative semigroup $S$, then $S$ has the group of fractions $G_S$. 

$\lhd$
By Proposition 2.1, we have $S = \Delta(S)$. Since $\Delta(S)$ is a right and left Ore set, $S$ has the group of fractions $G_S$.
$\rhd$

It follows from Corollary 2.2 that, in the study of centrally essential semigroup algebras of cancellative semigroups, it is sufficient to consider only semigroups $S$ which have the group of fractions $G_S$.

\textbf{Corollary 2.3.}\label{cor2.3}
Let $F$ be a field of characteristic $0$. Then any centrally essential semigroup algebra $FS$ of the cancellative semigroup $S$ over $F$ is commutative.

$\lhd$ The algebra $FS$ is semiprime if and only if is semiprime algebra $FG_S$; see \cite[Theorem 7.19]{Ok91}. It is well known that any group algebra over the field of characteristic $0$ is semiprime; e.g., see \cite[Theorem 4.2.12]{P77}. In \cite[Proposition 3.4]{MT18}, it is proved that centrally essential semiprime rings are commutative.~$\rhd$

\textbf{Example 2.1.}\label{ex2.1}
For a field $F$ of characteristic $0$ and the ring $M_{7}(F)$ of all matrices of order 7, we consider the subring $\mathcal{R}$ in $M_{7}(F)$ consisting of matrices of the form 
$$
\left(\begin{matrix}
\alpha & a & b & c & d & e & f\\
0 & \alpha & 0 & b & 0 & 0 & d\\
0 & 0 & \alpha & 0 & 0 & 0 & e\\
0 & 0 & 0 & \alpha & 0 & 0 & 0\\
0 & 0 & 0 & 0 & \alpha & 0 & a\\
0 & 0 & 0 & 0 & 0 & \alpha & b\\
0 & 0 & 0 & 0 & 0 & 0 & \alpha\\
\end{matrix}\right).
$$
Then $\mathcal{R}$ is a non-commutative centrally essential ring; see \cite[Example 2.4]{LT2la}. 
Let $e_{\alpha} = E_7$, $e_a$, $e_b$, $e_c$, $e_d$, $e_e$, $e_f$ be matrices, in which only non-null entry with value $1$ in places $\alpha$, $a$, $b$, $c$, $d$, $e$, $f$, respectively. We consider the semiring $S = \{e_{\alpha}, e_a, e_b, e_c, e_d, e_e, e_f\}\cup \{\theta\}$, where $\{\theta\}$ acts as zero. Then $\mathcal{R}\cong F_0S$, where $F_0S$ is the compressed semigroup algebra of the semigroup $S$ over the field $F$. Since 
$FS\cong F\bigoplus F_0S$ (e.g., see \cite[Corollary 4.9]{Ok91}), the  semigroup algebra $FS$ is centrally essential, since it is the direct sum of centrally essential algebras.

\subsection{The completion of the proof of Theorem 1.1}

\textbf{1.} Let $FS$ be a centrally essential ring and $0\neq a\in FG_S$, $a = \sum_{i=1}^n\alpha_ig_i$, where $\alpha_i\in F$, $g_i\in G_S$. It is known that $G_S = SZ(S)^{-1}$; see \cite[Proposition 9.8(iv)]{Ok91}. Then $a = \sum_{i=1}^n\alpha_is_it_i^{-1}$ for some $s_i\in S$, $t_i\in Z(S)$, $i = 1,\ldots, n$. We set $a' = \alpha_1s_1t_2\ldots t_n +\ldots + \alpha_ns_nt_1\ldots t_{n-1}\in FS$. We note that $a'\neq 0$. For this purpose, it is sufficient to verify that $s_1t_2\ldots t_n,\ldots, s_nt_1\ldots t_{n-1}$ are distinct elements in $FS$. Indeed, if $s_it_1\ldots\widehat{t_i}\ldots t_n = s_jt_1\ldots\widehat{t_j}\ldots t_n$ for $i\neq j$, then we multiply this equality by $(t_1\ldots t_n)^{-1}$ and obtain $s_it_i^{-1} = s_jt_j^{-1}$, i.e., $g_i = g_j$. This is a contradiction. By assumption, $0\neq a'c' = d'$ for some $c', d'\in Z(FS)$. Then $0\neq ac'' = d'$, where $c'' = t_1\ldots t_nc'$ and $d'$ are central elements in $FS$ which remain central in $FG_S$; see \cite[Corollary 9.11(i)]{Ok91}.

Conversely, let $0\neq a\in FS$, $a = \sum\alpha_is_i$, where $\alpha_i\in F$, $s_i\in S$. By assumption, $0\neq ac = d$ for some $c, d\in Z(FG_S)$,
$c = \sum_{i=1}^n\beta_ig_i$, $d = \sum_{j=1}^m\gamma_jh_j$, where $g_i, h_j\in G_S$. Let $g_i = x_iy_i^{-1}$, $h_j = z_jt_j^{-1}$ and $x_i, z_j\in S$, 
$y_i, t_j\in Z(S)$, $i = 1,\ldots, n$, $j = 1,\ldots, m$. We set 
$$
y = y_1\ldots y_n,\; t = t_1\ldots t_m,\;c' = cyt\in Z(FS).
$$
Then
$$
ac' = (ac)yt = dyt\in Z(FS).
$$ 
It remains to verify that $ac'\neq 0$. We have that 
$$
dyt = \gamma_1z_1yt_2\ldots t_m + \ldots + \gamma_mz_myt_1\ldots t_{m-1}.
$$
If $i\neq j$ and $z_iyt_1\ldots\widehat{t_i}\ldots t_m = z_jyt_1\ldots\widehat{t_j}\ldots t_m$, then 
$z_it_i^{-1} = z_jt_j^{-1}$ and $g_i = g_j$; this is a contradiction. 

\textbf{2.} The assertion follows from Example 2.1.~\hfill$\square$

\textbf{Example 2.2.}
Let $S = \langle x,y,z\mid \, z\in Z(S), z^2 = e, xy = zyx\rangle$. It is directly verified that $S$ is a cancellative semigroup which has the group of fractions $G_S = \langle x,y,z\mid \, z\in Z(G_S), z^2 = e, x^{-1}y^{-1}xy = z\rangle$. Since $z$ is a central involution and $x^2, y^2\in Z(G_S)$,  the unique non-trivial commutator in $G_S$ is $x^{-1}y^{-1}xy$. Therefore, the commutant $G'_S = <z>$. We have that $Z(G_S) = <x^2, y^2, z>$. Let 
$H = G'_S = \{e, z\}$, $\text{char }F = 2$ and $\hat{H} = e + z$. We verify that for $0\neq \alpha\in FG_S$, we have $\alpha\hat{H} = \beta\in Z(FG_S)$. Indeed, if $\alpha = {\underset{g\in \text{supp }(\alpha)}{\sum}}a_gg$, then $\beta = {\underset{g\in \text{supp }(\alpha)}{\sum}}a_gg\hat{H}$. Then for any $x\in G_S$ and $g\in \text{supp}(\alpha)$, we have
$$
[x, g\hat{H}] = [x, g]\hat{H} = xg(1 - g^{-1}x^{-1}gx)\hat{H} = 0,
$$
since $g^{-1}x^{-1}gx\in G'\subseteq Z(G_S)$.
If $\alpha\hat{H} = 0$, then $\alpha\in FG_S\hat{H}$; see \cite[Lemma 3.1.2]{P77}. In this case $\alpha\in Z(FG_S)$. Consequently, group algebra $FG_S$ is centrally essential. By Theorem 1.1, the semigroup algebra $FS$ is centrally essential as well.

\textbf{Example 2.3.}
Let $S = <x, y, z\mid \, z\in Z(S), xy = zyx>$. The semigroup $S$ has the group of fractions $G_S$ which is a free nilpotent group of nilpotence class 2; see \cite[Example 24.21]{Ok91}. It is known that if a group does not contain elements of order $p$, then centrally essential group algebra is commutative; see \cite[Proposition 1]{LT21b}. Therefore, the group algebra $FG_S$ is not centrally essential. By Theorem 1.1, the semigroup algebra $FS$ is not centrally essential as well.

\section{Proof of Theorem 1.2}\label{sec3}

\textbf{Lemma 3.1, \cite[Proposition 3]{LT21b}.}\label{lemma3.1}
Let $R$ be a ring. If for any non-zero-divisor $b$, there exists a non-zero-divisor $x$ such that $bx\in Z(R)$, then $R$ has the right classical ring of fractions.

\textbf{Lemma 3.2.}\label{lemma3.2} 
Let $FS$ be a centrally essential semigroup algebra of the cancellative semigroup $S$. Then for every non-zero-divisor $b\in FS$, there exists a non-zero-divisor $z\in FS$ such that $bz\in Z(FS)$.

$\lhd$ It follows from \cite[Lemma 4.4.4]{P77} that there exists a non-zero-divisor $x\in FG_S$ such that $bx = y\in Z(FG_S)$. If 
$x = \sum_{i=1}^n\alpha_is_it_i^{-1}$, where $\alpha_i\in F$, $s_i\in FS$, $t_i\in Z(FS)$, $i = 1,2, \ldots, n$, then element $z = xt_1\ldots t_n$ is a non-zero-divisor in $FS$ and $bz\in Z(FS)$.~$\rhd$

\textbf{Remark 3.3.}
The assertions of Lemma 3.1 and Lemma 3.2 are true for left classical rings of fractions as well. In this case, for a non-zero-divisor $b$, there exists a non-zero-divisor $x$ such that $xb$ is a central element. 

\textbf{Proposition 3.4.}\\ 
If $FS$ is a centrally essential semigroup algebra of the cancellative semigroup $S$, then $FS$ has the classical ring of fractions.

$\lhd$
The assertion follows from Lemma 3.1, Lemma 3.2 and Remark 3.3.~$\rhd$

\subsection{The completion of the proof of Theorem 1.2}

Let $FS$ be a centrally essential ring and $0\neq as^{-1}\in Q_{\text{cl}}(FS)$, where $s$ is a non-zero-divisor in $FS$. Let $\gamma\in FS$ be a non-zero-divisor such that $s\gamma = t\in Z(FS)$. Then $s^{-1} = \gamma t^{-1}$ in the ring $Q_{\text{cl}}(FS)$. By assumption, for the element $a\gamma\in FS$, there exist non-zero elements $c, d\in Z(FS)$ such that $0\neq (a\gamma)c = d\in Z(FS)$. Then
$$
(as^{-1})c = (a\gamma t^{-1})c = (a\gamma c)t^{-1} = dt^{-1}\neq 0,
$$
where $dt^{-1}\in Z(Q_{\text{cl}}(FS))$. Consequently, $Q_{\text{cl}}(FS)$ is a centrally essential ring.

Conversely, let $0\neq s\in FS$. By assumption, there exist two elements $t, r\in Z(Q_{\text{cl}}(FS))$ such that $0\neq st = r$. We note that
$Z(Q_{\text{cl}}(FS))\subseteq Q_{\text{cl}}(Z(FS))$; cf., \cite[Theorem 4.4.5]{P77}. Indeed, let $\rho\in Z(Q_{\text{cl}}(FS))$, 
$\rho = \alpha\beta^{-1}$, where $\alpha, \beta\in FS$ and $\beta$ is a non-zero-divisor. Then $\alpha\beta = \beta\alpha$ and $\alpha\beta^{-1} = \beta^{-1}\alpha$. By Lemma 3.2, there exists a non-zero-divisor $\gamma\in FS$ such that $\beta\gamma\in Z(FS)$. We set 
$\varepsilon = \beta\gamma$, $\eta = \alpha\gamma$ and obtain
$$
\eta\epsilon^{-1} = \alpha\gamma\gamma^{-1}\beta^{-1} = \alpha\beta^{-1} = \rho.
$$
In addition, $\varepsilon, \eta\in Z(Q_{\text{cl}}(FS))$. Considering this, we get $t = cd^{-1}$, $r = mn^{-1}$ for some $c, d, m, n\in Z(FS)$. Then
$$
s(cn) = (sc)n = (mn^{-1}d)n = md\in Z(FS),
$$
and $md\neq 0$, since $d$ is a non-zero-divisor in $FS$.~\hfill$\square$


\begin{thebibliography}{999}
\bibitem{CP61} Clifford, A.H., Prieston, G.B.: The Algebraic Theory of Semigroups. AMS Survey No. 7, Providence,  Volume 1 (1961)

\bibitem{LT20} Lyubimtsev, O.V., Tuganbaev, A.A.: Centrally essential endomorphism rings of abelian groups. Comm. Algebra. {\bf 48}(3),  
1249--1256 (2020)

\bibitem{LT2la} Lyubimtsev, O.V., Tuganbaev, A.A.: Centrally Essential Torsion-Free Rings of Finite Rank. Beitr\"age zur Algebra und Geometrie  Contributions to Algebra and Geometry. {\bf 62}(3),  615--622 (2021)

\bibitem{LT21b} Lyubimtsev, O.V., Tuganbaev, A.A.: Centrally essential group algebras and classical rings of fractions. Lobachevskii Journal of Mathematics. {\bf 42}(12), 2890--2894 (2021)

\bibitem{MT18} Markov, V.T., Tuganbaev, A.A.:  Centrally essential group algebras. J.~Algebra. {\bf 512}(15), 109--118 (2018)

\bibitem{MT19b} Markov, V.T., Tuganbaev, A.A.: Rings essential over their centers. Comm. Algebra. {\bf 47}(4), 1642--1649 (2019)

\bibitem{MT20b} Markov, V.T., Tuganbaev, A.A.: Uniserial Noetherian Centrally Essential Rings. Comm. Algebra. {\bf 48}(1), 149--153  (2020)

\bibitem{Ok91} Okninski, J.: Semigroup Algebras. Dekker, New York and Basel (1991)

\bibitem{P77} Passman, D.S.: The Algebraic Structure of Group Rings. John Wiley and Sons, New York (1977)
\end{thebibliography}
\end{document}